# Numerical Simulation of Dendritic crystal growth using phase field method and investigating the effects of different physical parameter on the growth of the dendrite.

Rahul Sanal
*Department of Mechanical and Aerospace Engineering, University at Buffalo, SUNY*
*rahulsan@buffalo.edu*

*Abstract*— **In this article, we study the phase-field model of solidification for numerical simulation of dendritic crystal growth that occurs during the casting of metals and alloys based on the kobayashi [1] model. Qualitative relationships between shapes of the crystal and physical parameters are studied and visualized.**

## I. INTRODUCTION

Solidification and dendritic growth is very important from a practical point of view. Some properties of solids, e.g. ductility, electrical conductivity and mechanical strength, are determined by the microscopic structures produced upon solidification. One would like to gain control of the structure formation to obtain the desired properties in manufacture. Applications can be found in for example casting and semiconductor production.

Dendrites are formed when surface anisotropy is included in the system, which means that there are going to be some preferred directions for solidification. Though there are many models to simulate dendritic growth in the literature but in the present report, a phase field model suggested by Ryo Kobayashi [1] is being studied.

The foundation of most solidification theories is the time-dependent Stefan problem. This theory describes the evolution of the thermal diffusion field around the solidification front with two conditions at the solid-liquid interface. The first condition relates the velocity of the moving front to the difference in thermal fluxes across the solid-liquid interface. The second, called Gibbs-Thomson condition, relates the interfacial temperature to the thermodynamic equilibrium, the local interfacial curvature and interface kinetics. Solving Stefan problem has mathematical difficulty of specifying conditions at an interface whose location is part of the unknown solution [2].

The phase-field method introduces an auxiliary continuous order parameter, $\varphi$ which is a function of position and time, to describe whether the material is liquid ($\varphi=0$) or solid ($\varphi=1$). The behavior of this variable is governed by an equation that couples to the evolution of the thermal field. In phase-field calculations, the conditions at the interface are not required. The location of the interface is obtained from the numerical solution for the phase-field variable at positions where $\varphi$ between 0 and 1. The phase-field method is attractive for numerical simulation since it allows one to solve a free-boundary problem without tracking the location of the interface.

For many materials, including metals, the surface energy of the solid-liquid interface and the kinetic coefficient depend on orientation. Since anisotropy has a crucial impact on the shape of the dendrites it is necessary to modify the phase-field model for simulation in realistic physical settings.

## II. MATHEMATICAL MODEL

The model includes two variables; one is a phase field $\phi(r, t)$ and the other is a temperature field $T(r, t)$. The variable $\phi(r, t)$ is an ordering parameter at the position r and the time t, $\phi = 0$ means being liquid and $\phi = 1$ solid. And the solid/liquid interface is expressed by the steep layer of $\phi$ connecting the values 0 and 1. Fig. 1 shows how the shape of crystal is described by the phase field $\phi$. In order to keep the profile of $\phi$ such form and to move it reasonably, we consider the following Ginzburg-Landau type free energy functional F similar to equation (2) including m as a parameter:

$$F = \int v \, [f(\phi, m) + (\epsilon^2 \phi/2)*|\nabla\phi|^2] \, dv \qquad (1)$$

where $\epsilon$ is a small parameter which determines the thickness of the layer. It is a microscopic interaction length and it also controls the mobility of the interface. f is a double-well potential which has local minimums at $\phi = 0$ and 1 for each m. Here we take the specific form of f as follows:

$$f(\phi, m) = 1/4\phi^4 - (1/2 - 1/3*m)\phi^3 + (1/4 - 1/2*m)\phi^2 \qquad (2)$$

Anisotropy can be introduced by assuming that $\epsilon$ depends on the direction of the outer normal vector at the interface. So $\epsilon$ is represented as a function of the vector $v = vi$ satisfying $\epsilon(\lambda,v) = \epsilon(v)$ for $\lambda>0$. The outer normal vector is represented by $-\nabla\phi$ at the interface. Thus, we consider:

$$F = \int v \, [f(\phi, m) + (\epsilon(-\nabla\phi)2/2)*|\nabla\phi|2] \, dv \qquad (3)$$

From the formula $\tau \, \partial\phi/\partial t = \partial f/\phi$ and further simplifying, we have the following evolution equation:

$$\tau \, \partial\phi/\partial t = -\partial/\partial x(\epsilon\epsilon' \partial\phi/\partial y) + \partial/\partial y(\epsilon\epsilon' \partial\phi/\partial x) + \nabla.(\epsilon^2 \nabla\phi) + \phi(1-\phi)(\phi-0.5+m) \qquad (4)$$

where $\tau$ is a small positive constant and $\partial\epsilon/\partial v=(\partial\epsilon/\partial vi)i$. The parameter m gives a thermodynamical driving force.

Especially in two dimensional space, we can take $\varepsilon = \varepsilon(\theta)$ where $\theta$ is an angle between v and a certain direction (for example the positive direction of the x-axis). $\varepsilon'$ means derivative with respect to $\theta$. Equation (4) gives the evolution of the order parameter or phase field variable $\phi$ with time.

Here we assume that m is a function of the temperature T, for example, $m(T) = \gamma(T_e - T)$ where $T_e$ is an equilibrium temperature, which means that the driving force of interfacial motion is proportional to the supercooling there. But in the following simulations, we used the form $m(T)=(\alpha-l[\gamma(T_e - T)]$ where $\alpha$ and $\gamma$ are positive constants; $\alpha < 1$, since this assures $|m(T)| < \frac{1}{2}$ for all values of T. Also m(T) is almost linear for T near $T_e$. To take anisotropy into account, let us specify $\varepsilon$ to be:

$$\varepsilon = 1 + \delta\cos(j(\theta-\theta_o)) \qquad (5)$$

The parameter $\delta$ means the strength of anisotropy and j is a mode number of anisotropy. Side branching can be stimulated in the dendrites by adding a small random noise in the equation (20). The noise can be of the form $a\phi(1-\phi)x$, where a is the strength of noise and x is a random number in the range[-0.5,0.5].

The equation for T is derived from the conservation law of enthalpy as:

$$\partial T/\partial t = \nabla^2 T + K\partial p/\partial t \qquad (6)$$

T is non-dimensionalized so that the characteristic cooling temperature is 0 and the equilibrium temperature is 1. K is a dimensionless latent heat which is proportional to the latent heat and inversely proportional to the strength of the cooling. For simplicity, the diffusion constant is set to be identical in both of solid and liquid regions (5) Is a heat conduction equation having a heat source along the moving interface, since $K\partial p/\partial t$ has non-zero value only when the interface passes through the point.

## III. SIMULATIONS

The simplest finite difference scheme with a nine point laplacian is used to solve equations (4) and (6). First, the new value of phase field variable is calculated to substitute it into the temperature field to get the new value of temperature field. Nine point laplacian is required for the stability of the solution.

A code in Matlab was developed to evolve the phase field variable and temperature field using equations (4) and (6) respectively and periodic boundary conditions.

The inputs needed for the simulation are as follows: nx, ny - size of the mesh

dx, dy - distance between the nodes in x & y direction dt - length of time step

timesteps - total number of timesteps

p(N,M) - initial phase field variable information T(N,M) - initial temperature field information k - latent Heat
tau - phase field relaxation time
epsilonbar - interfacial width
delta - strength of anisotropy
aniostropy - mode number of
anisotropy alpha - positive constant
gamma - positive constant

In all simulations the parameters are taken as follows: anisotropy= 4, on square mesh of size 500*500 with dt = 0.0001, latent heat k = 1.8, relaxation time $\tau$ = 0.0003, orientation angle $\theta$ = 1.57, interfacial width $\delta$ = 0.01 and dx = dy = 0.03 at time t = 0.2.

As we see in all the figures nucleation occurs at the center of the grid and it triggers the growth process.

The amplitude of stochastic noise is specified by $\alpha$ in the model. The side-branching propagation is known mainly due to thermal fluctuations, which enters solidification models in the form of random noise possessing specific features. Without noise, morphological instabilities at interfaces could never be expressed. Such noise stimulates fluctuations at the interface, which, when amplified, give rise to much of the structure observed in real systems. In other words, without noise, growth of the solid would be featureless, without the highly branched structures familiar in snowflakes. Therefore, to obtain more realistic simulation, $\alpha$ needs to be set to a suitable value.

## IV. EFFECT OF INTERFACIAL WIDTH ON DENDRITE GROWTH

In this subsection, for the same values of step size, *t* and a, we computed the problem for two different $\delta$ = 0.01 and 0.011 for time *t* = 0.12 as shown in Figure 1. We observed that when the interfacial width $\delta$ is slightly increased from 0.01 to 0.011 respectively, the tip velocity and tip radius of side branches also increased. Note that for interfacial width $\delta$ = 0 the growth of the crystal almost ceases.

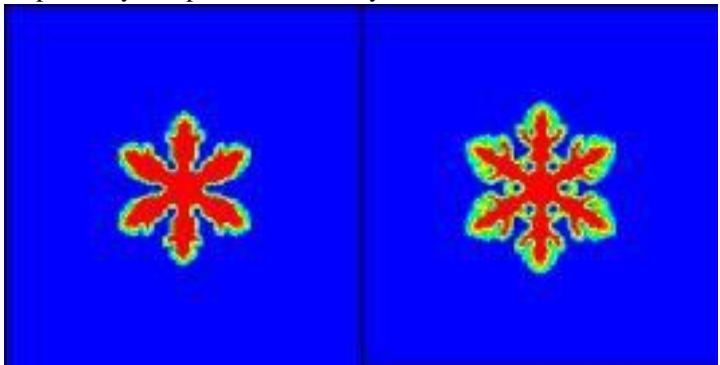

Figure 1. Dendtritic growth at $\delta$=0.01 and 0.011

## V. Anisotropy effects on Dendrite growth

Anisotropy has a direct effect on the shape of dendritic crystal. In Figure 2, we presented the tetragonal shape of dendrite using mode number of anisotropy= 4, on square mesh of size 500*500 with dt = 0.0001, latent heat k = 1.8, relaxation time τ = 0.0003, orientation angle θ = 1.57, interfacial width δ = 0.01 and dx = dy = 0.03 at time t = 0.2 both for phase-field and temperature values respectively. We now present the hexagonal shape of dendrite in Figure 2, by increasing the mode number of anisotropy= 6, both for phase-field and temperature values respectively at time t = 0.2. Note that anisotropy has a direct effect on the side branches. As clear from Figures 2, when anisotropy is increased from 4 to 6, the side branches appear around the crystal. Moreover the tip velocity and tip radius of side branches is also increased when anisotropy is changed from 4 to 6.

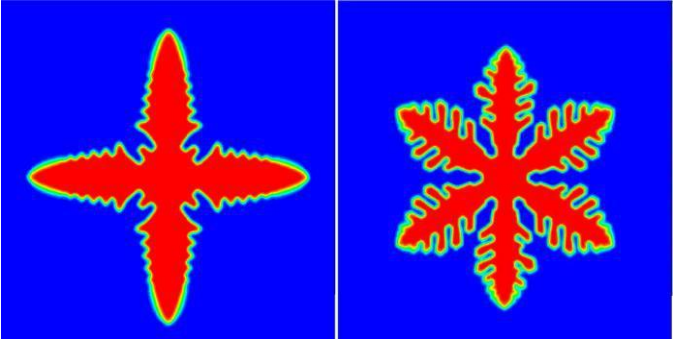

Figure 2. Dendtritic growth at anisotropy =4 and anisotropy =6

## VI. Latent Heat effects on Dendrite growth

The whole growth process of the dendrite depends mainly upon latent heat K. This auxiliary parameter is the most important physical quantity of the phase-field model. It accounts for the liberation of heat from the solidification front (diffusion layer). Large values of latent heat F accounts for the evacuation of heat from the interfacial region in a great amount. Therefore, the growth of dendrite is directly proportional to the value of latent heat F. For small values of F the growth process is slow, as less amount of heat evacuates from the diffusion layer and hence the solidification process is slow. Latent heat F is the main parameter which triggers the solidification process. Figures 3 shows the dendritic growth for latent heat F = 0.8, 1.0, 1.4, 1.8 and 2.0, inside square domain of size 500*500 with dt = 0.0002 , relaxation time τ = 0.0003 , orientation angle θ = 1.57 , interfacial width δ = 0.01 , mode number of anisotropy a = 6 and dx = dy = 0.03 respectively at time t = 0.1 . Figures 3 clearly shows that the shape of dendrite depends upon the latent heat k, for large values of k the branches and side branches form around the dendritic crystal.

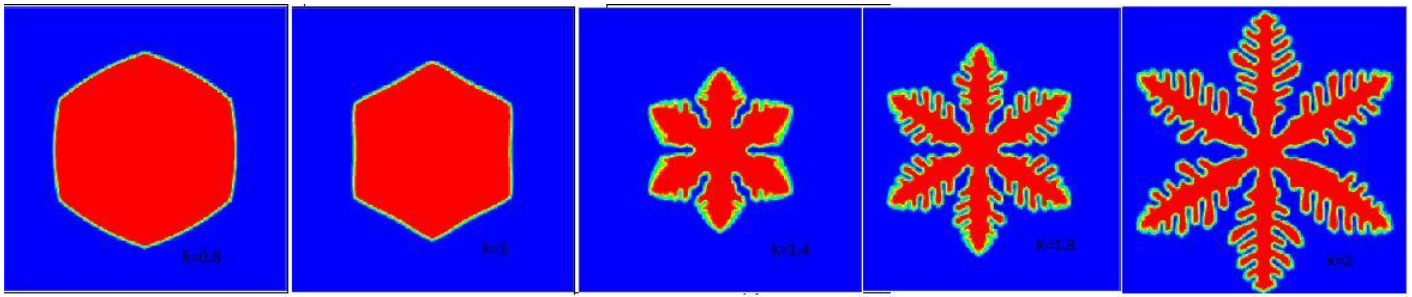

Figure 3 Dendritic growth with different latent heat values

## VII. Conclusion

In this work, a model is developed for the solution of phase-field equation coupled with thermal diffusion equation. The phase field method is analyzed by studying the effect of different physical parameter on the growth of dendrite is analyzed and justified with the actual phenomenon of dendritic solidification.

## VIII. APPENDIX

```matlab
clear all; clc;

%========================================
%========================================
 nx=500;
 ny=500;
dx=0.03;                    %e x direction
dy=0.03 ;                   %grid size y direction
dt=0.0003 ;                 %time step
tau=0.0003;
epsilonbar= 0.01;
mu=1.0;
k=4  ;
delta=0.02;
anisotropy= 4.0 ;
alpha= 0.9;
gamma= 10.0;
teq=1.0 ;

x= 1:nx;
y= 1:ny;
t=zeros(nx,ny) ; epsilon = zeros(nx,ny);
epsilon_derivative = zeros(nx,ny); phi1= zeros(nx,ny);
grad_phi_x=zeros(nx,ny);
grad_phi_y=zeros(nx,ny);
lap_phi=zeros(nx,ny);
lap_t=zeros(nx,ny);
phinew=zeros(nx,ny);
tnew=  zeros(nx,ny);
angl=  zeros(nx,ny);
times=0;

%====================setting the nuclei==================================
for i=1:nx
for j=1:ny

phi1(i,j)=0;

if ((i-nx/2)*(i-nx/2)+(j-ny/2)*(j-ny/2)<20);
phi1(i,j) = 1; end

end
end
phi=phi1;

%=====================gradient&Laplace===================================

for times=1:2000000

for i=1:nx
for j=1:ny

%periodic boundary condition jp = j+1;
jm = j-1;
ip = i+1;
im = i-1;
   if (im==0)
      im=nx;
   elseif (ip==nx+1)
      ip=1;
   end
   if (jm==0)
      jm=ny;
   elseif (jp==ny+1)
      jp=1;
   end

%gradient
grad_phi_x(i,j) = (phi(ip,j) - phi(im,j))/dx;
grad_phi_y(i,j) = (phi(i,jp) - phi(i,jm))/dy;

%laplacian
lap_phi(i,j) = (2.0*(phi(ip,j)+phi(im,j)+phi(i,jp)+phi(i,jm)) + phi(ip,jp)+phi(im,jm)+phi(im,jp)+phi(ip,jm) - 12.0*phi(i,j))/(3.0*dx*dx);
lap_t(i,j) = (2.0*(t(ip,j)+t(im,j)+t(i,jp)+t(i,jm)) + t(ip,jp)+t(im,jm)+t(im,jp)+t(ip,jm) - 12.0*t(i,j))/(3.0*dx*dx);

   if (grad_phi_x(i,j)==0) if (grad_phi_y(i,j)<0 )
       angl(i,j) =-0.5*pi;
      elseif (grad_phi_y(i,j)>0 )
       angl(i,j) = 0.5*pi;
     end
   end

   if (grad_phi_x(i,j)>0) if (grad_phi_y(i,j)<0)
       angl(i,j)= 2.0*pi + atan(grad_phi_y(i,j)/grad_phi_x(i,j));
      elseif (grad_phi_y(i,j)>0)
       angl(i,j)= atan(grad_phi_y(i,j)/grad_phi_x(i,j));
     end
   end
   if (grad_phi_x(i,j)<0)
       angl(i,j) = pi + atan(grad_phi_y(i,j)/grad_phi_x(i,j));
     end
     epsilon(i,j) = epsilonbar*(1.0 + delta*cos(anisotropy*(angl(i,j))));
     epsilon_derivative(i,j) = -epsilonbar*anisotropy*delta*sin(anisotropy*(angl(i,j)));
     grad_epsilon2_x = (epsilon(ip,j)^2 - epsilon(im,j)^2)/dx;
     grad_epsilon2_y = (epsilon(i,jp)^2 - epsilon(i,jm)^2)/dy;
```

```
        end
    end

%===================evolution
=======================================
=====

for i=1:nx
for j=1:ny

% periodic boundary
condition jp = j+1;
jm = j-1;
ip = i+1;
im = i-1;
 if (im==0)
    im=nx;
 elseif (ip==nx+1)
    ip=1;
 end
 if (jm==0)
    jm=ny;
 elseif (jp==ny+1)
    jp=1;
 end
term1=
(epsilon(i,jp)*epsilon_derivative(i,jp)*
grad_phi_x(i,jp) - epsilon(i,jm)*epsilon
_derivative(i,jm)*g rad_phi_x(i,jm))/dy;

term2= -
(epsilon(ip,j)*epsilon_derivative(ip,j)*
grad_phi_y(ip,j) - epsilon(im,j)*epsilon
_derivative(im,j)*g rad_phi_y(im,j))/dx;

term3= grad_epsilon2_x*grad_phi_x(i,j)+
grad_epsilon2_y*grad_phi_y(i,j);
```

```
phiold = phi(i,j);
m = alpha/pi * atan(gamma*(teq-t(i,j)));

%time evolution

phinew(i,j) = phi(i,j) + (term1 + term2
+ epsilon(i,j)*epsilon(i,j)*lap_phi(i,j)
+ term3          +          phiold*(1.0-
phiold)*(phiold-0.5+m))*dt/tau;
tnew(i,j) =t(i,j) + lap _t(i,j)*dt
+ k*(phi(i,j) - phiold);
phi(i,j)=phinew(i,j);
t(i,j)=tnew(i,j);

%visualization of the
output disp(phi);

figure(1);
image(phi*50);colormap('jet(64)');
pcolor(phi);shading flat; axis ('xy');

end
end
end
```